\documentclass[reqno, 12pt]{amsart}
\usepackage{amsmath, amssymb, amsthm, verbatim, bbm}

\newcommand \N {\mathbb{N}}
\newcommand \R {\mathbb{R}}
\newcommand \C {\mathbb{C}}
\newcommand \Z {\mathbb{Z}}
\newcommand \Sc {\mathcal{S}}

\newcommand \Oh {\mathcal{O}}
\newcommand \la {\langle}
\newcommand \ra {\rangle}
\newcommand \wt {\widetilde}
\newcommand \D {\partial}
\newcommand \eps {\varepsilon}

\newcommand{\ol}{\overline}
\newcommand \V {\mathcal{V}}
\DeclareMathOperator \re {Re}
\DeclareMathOperator \im {Im}

\DeclareMathOperator \sech {sech}
\newtheorem{prop}{Proposition}
\newtheorem{lem}{Lemma}
\newtheorem{thm}{Theorem}
\numberwithin{equation}{section}
\numberwithin{prop}{section}
\numberwithin{lem}{section}

\title
[Solitary waves for the Hartree equation with a slowly varying potential]
{Solitary waves for the Hartree equation with a slowly varying potential}

\author[K. Datchev]
{Kiril Datchev}
\email{datchev@math.berkeley.edu}
\author[I. Ventura]
{Ivan Ventura}
\email{iventura@math.berkeley.edu}
\address{Mathematics Department, University of California \\
Evans Hall, Berkeley, CA 94720, USA}

\begin{document}

\begin{abstract}We study the Hartree equation with a slowly varying smooth potential, $V(x) = W(hx)$, and with an initial condition which is $\eps \le \sqrt h$ away in $H^1$ from a soliton. We show that up to time $|\log h|/h$ and errors of size $\eps + h^2$ in $H^1$, the solution is a soliton evolving according to the classical dynamics of a natural effective Hamiltonian. This result is based on methods of Holmer-Zworski, who prove a similar theorem for the Gross-Pitaevskii equation, and on spectral estimates for the linearized Hartree operator recently obtained by Lenzmann. We also provide an extension of the result of Holmer-Zworski to more general inital conditions. \end{abstract}

\maketitle
\vspace{-1cm}

\section{Introduction}

In this paper we study the Hartree equation with an external potential:
\begin{equation}\label{e:hartree}
\begin{cases}i \D_t u = - \frac 12 \Delta u + V(x) u - \left(|x|^{-1} *|u|^2\right)u \\ u(x,0) = u_0(x) \in H^1(\R^3;\C).\end{cases}
\end{equation}
In the case $V \equiv 0$, solving the associated nonlinear
eigenvalue equation,
\begin{equation} \label{e:stationary} - \frac 12 \Delta \eta - \left(|\eta|^2 * \frac 1 {|x|}\right)\eta = - \lambda \eta, \end{equation}
gives solutions to \eqref{e:hartree} with evolution $u(t,x) =
e^{i\lambda t}\eta(x)$.  It is known that \eqref{e:stationary} has a
unique radial, positive solution $\eta \in H^1(\R^3)$ for a given
$\lambda > 0$; see \cite{lie} and \cite[Appendix A]{len}, as well as Appendix A below. For convenience of exposition in this paper we take $\lambda$ such that $\|\eta\|^2_{L^2} = 2$, but this is not essential. Using the
symmetries of \eqref{e:hartree}, we can construct from this $\eta$
the following family of \emph{soliton solutions} to
\eqref{e:hartree} in the case $V \equiv 0$:
\[
u(x,t) = e^{ix\cdot v}e^{i|v|^2t/2}e^{i\gamma}e^{i\lambda t}\mu^2\eta(\mu(x -a -
vt)), \quad (a,v,\gamma,\mu) \in \R^3 \times \R^3
\times\R\times\R_+.
\]
If $V \not\equiv 0$ but is slowly varying, there exist approximate soliton solutions in a sense made precise by the following theorem.

\begin{thm} \label{thm:main}
Let $V(x) = W(hx)$, where $W \in C^3(\R^3;\R)$ is bounded together with all derivatives up to order 3. Fix a constant $0 < c_1$, and fix $(v_0,a_0)
\in \R^3 \times \R^3$. Suppose $0 < \delta \le 1/2$, $0 < h \le h_0$, and $u_0 \in H^1(\R^3)$
satisfies
\[
\|u_0 - e^{iv_0\cdot (x-a_0)}\eta(x-a_0)\|_{H^1}\le c_1 h^2.
\]
Then if $u(t,x)$ solves \eqref{e:hartree} and
\[0 \le t \le  \frac {c_1} h + \frac{\delta|\log h|}{c_2 h},\]
we have
\[
\left\|u(t,x) - e^{v(t)\cdot (x - a(t))}
e^{i\gamma(t)}\eta\big[(x-a(t))\big]\right\|_{H^1_x(\R^3)} \le c_2 h^{2-\delta}.
\]
Here $(a,v,\gamma)$ solve the following system of equations
\begin{equation}\label{odes} \dot a = v, \quad \dot v = -\frac 12 \int \nabla V\left(x+a\right) \eta^2 (x)dx,\end{equation}
\[  \dot \gamma = \frac 12 |v|^2 + \lambda - \frac 12 \int V\left(x + a\right) \eta^2(x)dx + \frac 12 \int x \cdot \nabla V\left(x + a\right) \eta^2(x)dx, \]
with initial data $(a_0,v_0,0)$. The constants $h_0$ and $c_2$, depend only on $c_1$, $|v_0|$, and $\|W\|_{C^3(\R^3)}$. They are in particular independent of $\delta$.
\end{thm}

Note that in \eqref{odes}, the equation of motion of the center of mass $a$ of the soliton is given by Newton's equation:
\[\ddot a = - \nabla \overline{V}(a),\]
where $\overline{V} \stackrel{\textrm{def}}= V * \eta^2/2$. Observe also that because $\eta$ is exponentially localized (see Appendix \ref{app:eta}), $\eta^2/2$ is an approximation of a delta function and hence the effective potential $\overline{V}$ which governs the motion of the soliton is an approximation of $V$. The more complicated evolution of $\gamma$ is explained by the Hamiltonian formulation of the problem developed in Section \ref{s:ham}.

Our next theorem gives a slightly weaker result in the case of a more general initial condition.

\begin{thm} \label{thm:main2}
Let $V(x) = W(hx)$, where $W \in C^3(\R^3;\R)$ is bounded together with all derivatives up to order 3. Fix constants $0 < c_1$, and $0 \le 2\delta \le \delta_0 < 3/4$, and fix $(v_0,a_0) \in \R^3 \times \R^3$. Suppose $ 0 < h \le h_0 $, and $u_0 \in H^1(\R^3)$
satisfies
\[
\|u_0 - e^{iv_0\cdot (x-a_0)}\eta(x-a_0)\|_{H^1} \stackrel{\textrm{def}}=\eps \le c_1 h^{\frac12 + \delta_0}.
\]
 Then for
\[0 \le t \le  \frac {c_1} h + \frac{\delta|\log h|}{c_2 h},\]
we have
\[
\left\|u(t,x) - e^{v(t)\cdot (x - a(t))}
e^{i\gamma(t)}\mu(t)^2\eta\big[\mu(t)(x-a(t))\big]\right\|_{H^1_x(\R^3)} \le c_2 h^{-\delta}\tilde\eps,
\]
where $\tilde \eps \stackrel{\textrm{def}}= \eps + h^2$. Here $(a,v,\mu,\gamma)$ solve the following system of equations
\[ \dot a = v + \Oh(\tilde\eps^2), \quad \dot v = -\frac\mu2 \int \nabla V\left(\frac x\mu+a\right) \eta^2 (x)dx + \Oh(\tilde\eps^2), \qquad \dot\mu = \Oh(\tilde\eps^2),\]
\[  \dot \gamma = \frac 12 |v|^2 + \lambda \mu^2 - \frac 12 \int V\left(\frac x \mu + a\right) \eta^2(x)dx - \frac 1{2\mu} \int x \cdot \nabla V\left(\frac x\mu + a\right) \eta^2(x)dx + \Oh(\tilde\eps^2), \]
with initial data $(a_0,v_0,1,0)$.
The constants $h_0$ and $c_2$, as well as the implicit constants in the $\Oh$ error terms, depend only on $c_1$, $|v_0|$, and $\|W\|_{C^3(\R^3)}$. They are in particular independent of $\delta$.\end{thm}

This phenomenon has been studied in the physics literature by Eboli-Marques \cite{em}, who show for various explicit (but not necessarily slowly varying) potentials $V$ that soliton solutions which obey Newtonian equations of motion exist. Similar theorems have been proven in the case of more general nonlinearities by Fr\"ohlich-Gustafson-Jonsson-Sigal \cite{fgjs} and by Fr\"ohlich-Tsai-Yau \cite{fty2}. More recently Jonsson-Fr\"ohlich-Gustafson-Sigal \cite{jfgs} have extended the validity of the effective dynamics to longer time in the case of a confining potential $V$, and Abou-Salem \cite{a} has treated the case of a potential $V$ which is permitted to vary in time. The case of a power nonlinearity was studied by Bronski-Jerrard \cite{bj}, and the case of the cubic nonlinear Schr\"odinger equation in dimension one was also studied by Holmer-Zworski \cite{hz1}, \cite{hz2}. Other papers have established effective classical dynamics in quantum equations of motion in a wide variety of settings: see \cite{fgjs} and \cite{a} for many references.

Our result improves those of \cite{fgjs} and \cite{a} in the case of the equation \eqref{e:hartree} in several respects.  First we provide a more precise error bound, improving $\tilde\eps$ from $h + \eps$ to $h^2 + \eps$. Second we remove the errors in the equations of motion in the case $\eps = \Oh(h^{2-\delta})$. Finally, we establish the effective dynamics for longer time: in \cite{fgjs} the result obtained was valid only up to time $c(\eps^2 + h)^{-1}$ for a small constant $c$, while in \cite{a} the result was valid only up to time $\delta |\log h|/h$ and required the assumption $\eps = \Oh(h)$.

In \cite{fgjs} more general initial data are considered, that is to say $\eps$ is assumed to be small but not necessarily $\Oh(h^{1/2 +})$, although in this case the result is obtained only for time $\eps^{-2}$. In that situation the methods of the present paper, although applicable, do not improve that result, so for ease of exposition we have considered only the special case $\eps = \Oh(h^{1/2+})$ where we have an improvement. 

In this paper we follow most closely \cite{hz2}, which in turn builds on \cite{hz1} and on earlier work on soliton stability going back to Weinstein \cite{w} (see those papers for more references). We adapt those arguments to a higher-dimensional setting where in particular there is no longer an explicit form for $\eta$, and to the nonlocal Hartree nonlinearity. For this last task we make use of the classical Hardy-Littlewood-Sobolev inequalty and of spectral estimates for the linearized Hartree operator
\[\mathcal{L} w \stackrel{\textrm{def}}= -\frac{1}{2} \Delta u - \left(\frac{1}{|x|} * \eta(w+\ol{w})\right)\eta - \left(\frac{1}{|x|} * \eta^2 \right)w+\lambda w,\]
obtained recently by Lenzmann \cite{len}.

We also extend the methods of \cite{hz2} in that we adapt them to more general initial data. It is at this point that our proofs depart most significantly from those of \cite{hz2}, and this work is contained in Section \ref{reparam}. The crucial additional element is a closer analysis of the differential equation for the error studied in Lemmas \ref{5.1} and \ref{5.2}. This closer analysis applies also to the Gross-Pitaevskii equation studied in \cite{hz2}, giving us Theorem \ref{thm:main3} below.

To state this theorem, we suppose $u\colon\R \times \R \to \C$ solves
\begin{equation}
\label{e:gp}
\begin{cases}
i \D_t u = - \frac 12 \D_x^2 u + V(x) u - |u|^2u, \\ u(x,0) = u_0(x) \in H^1(\R;\C).
\end{cases}
\end{equation}
In this case the ground state soliton solution of the corresponding elliptic nonlinear eigenvalue equation
\[-\frac 12 \eta = - \frac 12 \eta'' - \eta^3\] 
is given by
\[\eta(x) = \sech(x).\]
We then have
\begin{thm} \label{thm:main3}
Let $V(x) = W(hx)$, where $W \in C^3(\R;\R)$ is bounded together with all derivatives up to order 3. Fix constants $0 < c_1$, $0 < \delta_0 < 3/4$ and fix $(v_0,a_0) \in \R \times \R$. Suppose $0 \le 2\delta \le \delta_0$ and $0 < h \le h_0$. For $u_0 \in H^1(\R)$ put
\[\|u_0 - e^{iv_0\cdot (x-a_0)}\sech(x-a_0)\|_{H^1} \stackrel{\textrm{def}}= \eps \le c_1 h^{\frac 12 + \delta_0}\]
Then for
\[0 \le t \le  \frac {c_1} h + \frac{\delta|\log h|}{c_2 h},\]
we have
\[
\left\|u(t,x) - e^{v(t)\cdot (x - a(t))}
e^{i\gamma(t)}\mu(t)\sech\big[\mu(t)(x-a(t))\big]\right\|_{H^1_x(\R^3)} \le c_2 h^{-\delta}\tilde\eps,
\]
where $u$ solves \eqref{e:gp} and $\tilde \eps \stackrel{\textrm{def}}= \eps + h^2$. Here $(a,v,\mu,\gamma)$ solve the following system of equations
\[\dot a = v + \Oh(\tilde\eps^2), \quad \dot v = -\frac {\mu^2}2 \int V'(x + a)\sech^2(\mu x)dx + \Oh(\tilde\eps^2), \quad \dot\mu = \Oh(\tilde\eps^2),\]
\[\dot \gamma = \frac 12 \mu^2 + \frac 12 v^2 - \mu\int V \left(x + a\right)\sech^2(\mu x)dx + \mu^2\int xV\left(x + a\right)\sech^2(\mu x)\tanh(\mu x) dx + \Oh(\tilde\eps^2), \]
with initial data $(a_0,v_0,1,0)$.
The constants $h_0$ and $c_2$, as well as the implicit constants in the $\Oh$ error terms, depend only on $c_1$, $\delta_0$, $|v_0|$, and
$\|W\|_{C^3(\R^3)}$. They are in particular independent of $\delta$.
\end{thm}

To prove this result, one simply replaces Lemmas 5.1 and 5.2 of \cite{hz2} with Lemmas \ref{5.1} and \ref{5.2} of the present paper. Because the details are very similar to the ones given in Section \ref{reparam} below, we omit them.

The methods of this paper can be extended to the case of more general nonlinearities under additional spectral nondegeneracy assumptions: see \cite{fgjs} for examples. In that paper, and also in \cite{fty2}, more general classes of equations are considered under such assumptions. For the present work we have restricted our attention to two physical nonlinearities for which the necessary spectral results are known.

The outline of the proof and of this paper are as follows. 
\begin{itemize}
	\item In Section \ref{s:ham} we recast \eqref{e:hartree} as a Hamiltonian evolution equation in $H^1(\R^3)$, with the Hamiltonian given by \eqref{harham}. We define the manifold of solitons to be the set of functions of the form $e^{v\cdot(x-a)}e^{i\gamma}\mu^2\eta(\mu(x-a))$ for some $(a,v,\gamma,\mu) \in \R^3\times\R^3\times\R\times\R^+$, and we show that the equations \eqref{odes} come from the restriction of the Hamiltonian \eqref{harham} to this manifold.
	\item In Section \ref{spec} we review and extend slightly the relevant spectral results from \cite{len}.
	\item In Section \ref{reparam} we compute the differential equation for the difference between the true solution $u$ and the `closest point' on the manifold of solitons. We then estimate this difference, proving Thoerem \ref{thm:main2}.
	\item In Section \ref{final} we show how the additional assumption on the initial condition in Theorem 1 gives the exact equations of motion \eqref{odes}.
	\item Finally in Appendix \ref{app:eta} we collect the properties of $\eta$ which we need for our proofs, and in Appendix \ref{app:b} we review a standard proof of the global well-posedness of \eqref{e:hartree}.
\end{itemize}

\section{Hamiltonian equations of motion} \label{s:ham}

This section is divided into four subsections. In the first we define a symplectic structure on $H^1$ and recall a few basic lemmas from symplectic geometry. In the second we define the manifold of solitons, which has a natural action on it by the group of symmetries of \eqref{e:hartree}. We compute the Lie algebra associated to this group of symmetries and from that deduce a formula for the derivative of a curve in the group in terms of the Lie algebra. In the third we prove that the manifold of solitons is a symplectic submanifold and compute the restriction of the symplectic form to it. In the fourth we compute the Hartree Hamiltonian and its restriction to the manifold of solitons, and derive the equations \eqref{odes} as the equations of motion associated to the restricted Hamiltonian. Most of the ideas in this section are present in \cite[Section 2]{hz1}

\subsection{Symplectic Structure.} We work over the vector space
\[\mathcal{V} \stackrel{\textrm{def}}= H^1(\R,\C) \subset L^2(\R,\C),\]
viewed as a \emph{real} Hilbert space. The inner product and the
symplectic form are given by
\begin{equation}
\la u,v \ra \stackrel{\textrm{def}}= \re \int u \ol{v}, \qquad \omega
( u,v ) \stackrel{\textrm{def}}= \im \int u \ol{v},
\label{e:products}
\end{equation}

Let $H: \V \rightarrow \R$ be a function, a Hamiltonian. The
associated Hamiltonian vector field is a map $\Xi_H\colon\V \to T\V$.
The vector field $\Xi_H$ is defined by the relation
\begin{equation}
\omega(v,(\Xi_H)_u) = d_uH(v), \label{e:Hfield1}
\end{equation}

where $v \in T_u\V$, and $d_uH: T_u\V \rightarrow \R$ is defined by
$$
d_u H(v)=\frac{d}{ds} \Big|_{s=0} H(u+sv).
$$
In the notation above we have
\begin{equation}
d_uH(v) = \la dH_u,v \ra, \quad (\Xi_H)_u = -i dH_u,
\label{e:Hfield2}
\end{equation}
where the first equation provides a definition of $dH_u$, and the second a formula for computing $\Xi_H$.

For future reference present two simple lemmas from symplectic geometry.
The proofs for these can be found in \cite[Section 2]{hz1}.

\begin{lem} \label{lem:2.1}
Suppose that $g:\V \rightarrow \V$ is a diffeomorphism such that
$g^*\omega=\mu(g)\omega,$ where $\mu(g) \in C^\infty(\V,\R).$ Then
for $f \in C^\infty(\V,\R)$
\begin{equation}
(g^{-1})_* \left((\Xi_f)_{g(\rho)}\right) = \frac 1 {\mu(g)} \Xi_{g^*f}(\rho),  \qquad\rho \in \V.
\end{equation}
\end{lem}

Suppose that $f \in C^\infty(\V,\R)$ and that $df(\rho_0)=0$. Then
the Hessian of $f$ at $\rho_0$, $f''(\rho_0):T_\rho \V \mapsto
T^*_\rho \V$ is well defined. We can identify $T_\rho \V$ with
$T^*_\rho \V$ using the inner product, and define the Hamiltonian map
$F: T_\rho \V \rightarrow T_\rho \V$ by
\begin{equation} \label{e:hammap}
F= -i f''(\rho_0), \quad \la f''(\rho_0) X, Y \ra= \omega (Y ,F X).
\end{equation}
In this notation we have
\begin{lem} \label{lem:2.2}
Suppose that $N \subset V$ is a finite dimensional symplectic
submanifold of $V$ and $f \in C^\infty(V,\R)$ satisfies
\[
\Xi_f(\rho) \in T_\rho N \subset T_\rho V, \quad \rho \in N.
\]
If at $\rho_0 \in N$, $df(\rho_0) = 0$ then the Hamiltonian map
defined by \eqref{e:hammap} satisfies
\[
F(T_\rho N ) \subset T_\rho N.
\]
\end{lem}

\subsection{Manifold of solitons as an orbit of a group.}
For $g = (a,v,\gamma,\mu) \in \R^3\times\R^3\times\R\times\R_+$, we define the map
\begin{equation}\label{e:groupact}
H^1 \ni u \mapsto g \cdot u \in H^1, (g \cdot u)(x)
\stackrel{\textrm{def}}= e^{i\gamma} e^{iv(x-a)}\mu^2 u(\mu(x-a)).
\end{equation}
This action gives the following group structure on $\R^7 \times \R_+$:
\[(a,v,\gamma,\mu)\cdot(a',v',\gamma',\mu') = (a'',v'',\gamma'',\mu''),\]
where
\[
v'' = v + \mu v',\quad a'' = a + \frac{a'}\mu,\quad \gamma'' =
\gamma + \gamma' + \frac{va'}\mu,\quad \mu'' = \mu\mu'.
\]

The action of $G$ is conformally symplectic in the following sense:
\begin{equation}\label{confsymm}g^* \omega = \mu \omega,\quad g = (a,v,\gamma,\mu), \end{equation}
as is easily seen from \eqref{e:products}.

The Lie algebra of $G$, denoted $\mathfrak{g}$, is generated by the
following eight elements:

\begin{equation}\label{e:ei}\begin{array}{lll}
e_1 = -\D_{x_1},& e_4=ix_1 &  e_7 = i, \\
e_2 = -\D_{x_2},& e_5=ix_2, & e_8 = 2 + x \cdot \nabla. \\
e_3 = -\D_{x_3},& e_6=ix_3,\\
\end{array}\end{equation}

These are simply the partial derivatives at the identity of $(g\cdot
u)(x)$ with respect to each of the eight parameters
$(a,v,\gamma,\mu)$. The following computation gives the derivative
of a curve in $G$ in terms of this basis.

\begin{lem}\label{yeq}
Let $g \in C^1(\R, G)$ and $u \in \Sc(\R).$ Then, in the notation of
\eqref{e:groupact},
\[\frac d{dt} g(t)\cdot u = g(t)\cdot(Y(t)u),\]
where $Y(t) \in \mathfrak{g}$ is given by
\begin{equation}
\label{e:X} Y(t) = \mu(t) \sum_{j=1}^3\dot a_j(t)e_j +
\mu(t)\sum_{j=1}^3\frac{\dot v_j(t)}{\mu(t)}e_{3+j} + (\dot\gamma(t)
- \dot a(t)\cdot v(t))e_7 + \frac{\dot\mu(t)}{\mu(t)}e_8,
\end{equation}
where $g(t) = (a(t),v(t),\gamma(t),\mu(t)) =
(a_1(t),a_2(t),a_3(t),v_1(t),v_2(t),v_3(t),\gamma(t),\mu(t))$.
\end{lem}

We define the submanifold of solitons, $M \subset H^1$, as the orbit
of $\eta$ under $G$, where $\eta$ is the function described in Appendix \ref{app:eta}.
\begin{equation} \label{e:manifold}
M=G \cdot \eta \simeq G / \Z, \quad T_\eta M =\mathfrak{g} \cdot
\eta \simeq \mathfrak{g}.
\end{equation}
The quotient corresponds to the $\Z$-action
\[
(a,v,\gamma,\mu) \mapsto (a,v,\gamma+ 2 \pi k,\mu), \quad k \in \Z.
\]

We also record the following simple consequence of the implicit function theorem and of the nondegeneracy of $\omega$. The proof can be found, for example in \cite[Lemma 3.1]{hz1}.

\begin{lem}\label{3.1}
For $\Sigma$ and compact subset of $G/\Z$, let
\[
U_{\Sigma,\delta} = \{u \in H^1: \inf_{g \in \Sigma} \|u - g\cdot
\eta\|_{H^1} < \delta\}.
\]
If $\delta \le \delta_0 = \delta_0(\Sigma)$ then for any $u \in
U_{\Sigma, \delta}$, there exists a unique $g(u) \in \Sigma$ such
that
\[
\omega(g(u)^{-1} \cdot u - \eta, X \cdot \eta) = 0 \quad
\forall X \in \mathfrak{g}.
\]
Moreover, the map $u \mapsto g(u)$ is
in $C^1(U_{\Sigma, \delta}, \Sigma)$.
\end{lem}

\subsection{Symplectic structure on the manifold of solitons.}
We compute the symplectic form $\omega \big|_M$ on $T_\eta M$ by
using
$$
(\omega \big|_M)_\eta(e_i,e_j)=\im \int (e_i  \cdot \eta)(x)
(\ol{e_j \cdot \eta})(x).
$$
We take this opportunity to remind the reader (as mentioned in
Appendix \ref{app:eta}) that $\|\eta\|_{L^2}^2=2$ Using formulas
given in \eqref{e:ei} we compute all these forms.

\begin{lem}
The evaluation at $\eta$ of the restriction of the symplectic form to $M$ is given by
\[
\left(\omega \big|_M\right)_\eta =  (dv \wedge da + d\gamma \wedge
d\mu)_{(0,0,0,1)} =  (d(vda + \gamma d\mu))_{(0,0,0,1)}.
\]
\end{lem}

\begin{proof}
If $j,k$ are both taken from $\{1,2,3,8\}$ or both taken from
$\{4,5,6,7\}$, then the integrand $(e_j  \cdot \eta)(x) (\ol{e_k
\cdot \eta})(x)$ is a real function, implying that $(\omega
\big|_M)_\eta(e_j,e_k)=0$.

If $j \in \{1,2,3\}$ and $k \in \{4,5,6\}$ we have
$e_j=-\partial_j$ and $e_k=ix_{k-3}$.

\begin{itemize}
\item  If $j \neq k-3$ then integrating by parts gives
\[ (\omega \big|_M)_\eta(e_j,e_k)=\im \int (e_j  \cdot \eta)(x) (\ol{e_k \cdot \eta})(x) = \im \int (-\partial_j \eta)(\overline{ix_{k-3} \eta})  = -\int(\eta)(x_{k-3}  \partial_j \eta).\]
This implies that $(\omega \big|_M)_\eta(e_j,e_k)=0$

\item If $j=k-3$ by parts integration gives
\[ (\omega \big|_M)_\eta(e_j,e_k)=\im \int (e_j  \cdot \eta)(x)
(\ol{e_k \cdot \eta})(x)=\int (\partial_j \eta)(x_{j} \eta)=- \int (\eta(\eta+x_j \partial_j \eta)).
\]
Solving this yields that $(\omega \big|_M)_\eta(e_j,e_k) =-1$
\end{itemize}

If $j \in \{1,2,3\}$ and $k =7$ by parts integration gives
\[(\omega \big|_M)_\eta(e_j,e_k)=\im \int (e_j  \cdot \eta)(x)(\ol{e_k \cdot \eta})(x)=\im \int (-\partial_j \eta)(\overline{i \eta})= \int (\partial_j \eta)(\eta)= -\int(\eta)(  \partial_j \eta),\]
implying $(\omega \big|_M)_\eta(e_j,e_k)=0$.

If $j \in \{4,5,6\}$ and $k=8$, we get
\begin{eqnarray*}
(\omega \big|_M)_\eta(e_j,e_k)&=&\im \int (e_j  \cdot \eta)(x)
(\ol{e_k \cdot \eta})(x) =\im \int ix_j \eta (2+x \cdot \nabla) \eta \\
&=& 2\int x_j \eta^2 +\int x_j \eta x \cdot \nabla \eta \\
&=& 2\int x_j \eta^2  + \int x_j \eta (x_1 \partial_1 \eta +x_2
\partial_2 \eta + x_3 \partial_3 \eta).
\end{eqnarray*}
Now $\int x_j \eta^2$ is zero as it is odd in the $x_j$ variable. Since
all the terms in this last expression can be reduced to this for by integrating by parts
we see that $\left(\omega \big|_M\right)_\eta(e_j,e_k)=0$.

If $j=7$ and $k=8$ we observe that since by integration by parts we have $\int \eta
x \cdot \nabla \eta=-\frac{3}{2}\| \eta \|^2_{L^2}$, then
\[(\omega \big|_M)_\eta(e_j,e_k)=\im \int (e_j  \cdot \eta)(x)
(\ol{e_k \cdot \eta})(x)= \int \eta (2+ x \cdot \nabla) \eta= 2\|\eta\|^2_{L^2} - \frac{3}{2} \|\eta\|^2_{L^2},\]
giving that $(\omega \big|_M)_\eta(e_j,e_k)=1$.

Putting all this together gives the result.
\end{proof}

We now observe from \eqref{e:manifold} and \eqref{confsymm} that
\begin{equation} \label{e:symprest}
\omega \big|_M = \mu dv \wedge da + v d \mu \wedge da + d \gamma
\wedge d \mu.
\end{equation}

Now let $f$ be a function defined on $M$, $f=f(a,v,\gamma,\mu)$. The
associated Hamiltonian vector field, $\Xi_f$, is given by
\[
\omega(\cdot,\Xi_f)=df=f_a da + f_v dv+f_\mu d\mu+f_\gamma d \gamma.
\]

Using \eqref{e:symprest} we obtain
\begin{equation} \label{e:2.15}
\Xi_f = \frac{1}{\mu} \nabla_v f \cdot \nabla_a + \frac 1\mu
\left(-\nabla_a f - (\partial_\gamma f) v\right)\cdot \nabla_v +
\frac{\partial}{\partial \gamma}f
\partial_\mu + \left( \frac{1}{\mu} v \cdot \nabla_v f -
- \partial_\mu f \right)
\partial_\gamma.
\end{equation}
The Hamiltonian flow is obtained by solving
\[
\dot{v} =-\nabla_a f - (\partial_\gamma f) v, \quad \dot a = \frac 1
\mu \nabla_v f, \quad \dot \mu = \partial_\gamma f, \quad \dot
\gamma = \frac 1 \mu v \cdot \nabla_v f - \partial_\mu f.
\]

\subsection{The Hartree Hamiltonian restricted to the manifold of
solitons}
Using the symplectic form given in \eqref{e:products}, and
\[
H(u)\stackrel{\textrm{def}}=\int \frac{1}{4}|\nabla u|^2 - \frac{1}{4} |u|^2 \left(|u|^2 *
\frac{1}{|x|}\right),
\]
we find that
\[d_uH(v)= \re \int \left( -\frac{1}{2} \Delta u -\left(|u|^2 * \frac{1}{|x|}\right)u \right)\ol{v}.\]
The Hamiltonian flow associated to this vector field is
\begin{equation} \label{e:2.4}
\dot{u}=(\Xi_H)_u=-i \left( -\frac{1}{2} \Delta u - \left(|u|^2 * \frac{1}{|x|}\right)u \right).
\end{equation}

The restriction of
\[
H(u)=\int \frac{1}{4}|\nabla u|^2 - \frac{1}{4} |u|^2 \left(|u|^2 * \frac{1}{|x|}\right),
\]
to $M$ is given by computing
\[ H(g \cdot \eta) = \frac{|v|^2 \mu}{4} \|\eta \|^2_{L^2} + \mu^3 H(\eta) = \frac{|v|^2 \mu}{2} + \mu^3 H(\eta),\]
for $g = (a,v,\gamma,\mu)$. The flow of \eqref{e:2.15} for this $f$
describes the evolution of a soliton. We have in particular
\[\dot\gamma = \frac 12 |v|^2 - 3 \mu^2 H(\eta),\]
and because we know that $e^{i\lambda t} \eta(x)$ solves \eqref{e:hartree}, we can compute that $H(\eta) = -\lambda/3$.

We now consider the Hartree Hamiltonian,
\begin{equation}\label{harham}
H_V(u)=\frac{1}{4}\int |\nabla u|^2 - \frac{1}{4} \int |u|^2
\left(|u|^2 * \frac{1}{|x|}\right) + \frac 12 \int  V(x) |u|^2,
\end{equation}
and its restriction to $M = G \cdot \eta$ given by
\begin{equation} \label{e:mysterious}
H_V|_M = \frac{|v|^2 \mu}{2} + \lambda \frac{\mu^3}3 + \frac{\mu^4}{2}
\int V(x) \eta^2(\mu(x-a)).
\end{equation}

The flow of $H_V|_M$ can be read off from \eqref{e:2.15}:
\[
\dot{v} = -\frac\mu2 \int \nabla V\left(\frac x\mu+a\right) \eta^2 (x)dx, \qquad \dot a = v, \qquad \dot \mu =0,
\]
\[
\dot \gamma =  \frac 12 |v|^2 + \lambda \mu^2 - \frac 12 \int V\left(\frac x \mu + a\right) \eta^2(x)dx + \frac 1{2\mu} \int x \cdot \nabla V\left(\frac x\mu + a\right) \eta^2(x)dx.
\]
These are the same as the ones given in \eqref{odes}. The
evolution of $a$ and $v$ is simply the Hamiltonian evolution of
$\frac 12 |v|^2 + \frac{\mu^3}{2} \int \nabla V(x+a) \eta^2(\mu x) $ when
$\mu$ is held constant. As a result the evolution of the phase is
explained by \eqref{e:mysterious}.

Finally we give an important application of  Lemma \ref{lem:2.2}. We put
\[H_\lambda(u) = \int \frac{1}{4}|\nabla u|^2 - \frac{1}{4} |u|^2 \left(|u|^2 * \frac{1}{|x|}\right) + \frac \lambda 2 \int |u|^2,\]
and observe that $\eta$ is a critical point of this functional, while the Hessian of $H_\lambda $at $\eta$ is given by
\begin{equation}\label{el}\mathcal{L} w \stackrel{\textrm{def}}= -\frac{1}{2} \Delta u - \left(\frac{1}{|x|} * \eta(w+\ol{w})\right)\eta - \left(\frac{1}{|x|} * \eta^2 \right)w+\lambda w.\end{equation}
Now if in Lemma \ref{lem:2.2} we take, $H_\lambda$ to be $f$, $N$ to be the eight dimensional manifold of solitons $M$, and $\rho=\eta$, we find that
\begin{equation}\label{il} i\mathcal{L} \left(T_\eta M\right) \subset T_\eta M. \end{equation}

\section{Spectral estimates}\label{spec}

In this section we recall crucial spectral estimates for the
operator $\mathcal{L}$ from \eqref{el}, which is the linearization of $- \frac 12
\Delta u - \left(|u|^2 * \frac 1 {|x|}\right)u + \lambda u$. We observe that this operator can be decomposed as follows:
$$
\mathcal{L}w=\left[ \begin{array}{cc}
L_+&0\\0&L_-\end{array}\right]\left[ \begin{array}{c} \re w
\\ \im w \end{array}\right],
$$
with
$$
L_+ \re w=-\frac{1}{2} \Delta \re w- 2\left(\frac{1}{|x|}*\eta \re w
\right)\eta-\left(\frac{1}{|x|}*\eta^2 \right)\re w+ \lambda \re
w,
$$
and
$$
L_- \im w=-\frac{1}{2} \Delta \im w- \left(\frac{1}{|x|}*\eta^2
\right)\im w+\lambda \im w.
$$

From Remark 2 following Theorem 4 in \cite{len} we have the
following proposition:
\begin{prop}
Let $w \in H^1(\R,\C)$ and suppose that for any $X \in
\mathfrak{g}$, $\omega(w,X \eta)=0$. Then,
\begin{eqnarray}
\la \mathcal{L}w, w \ra \geq c \| w \|^2_{H^1}, \label{e:4.1}
\end{eqnarray}
where $c$ is an absolute constant.
\end{prop}
Now we consider solutions $f$ of the equation
\begin{equation}
L_+f=Q(x)\eta(x), \label{e:4.2}
\end{equation}
where $Q(x)$ is real-valued and of the form $Q(x)=a_0(t)+ \sum
a_{ij}(t)x_ix_j$, with $Q(x)\eta$ symplectically orthogonal to the
generalized kernel of $i \mathcal{L}$, and with $a_{ij}(t)$ bounded
in $t$.
\begin{prop}\label{4.2}
The equation \eqref{e:4.2} has a unique solution in $(\ker
(L^+))^\perp \subset L^2(\R^3)$. This solution is also in
$C^\infty(\R^3)$ with the property
\begin{equation}
\label{expdecayf}e^{\frac{1}{2}(\sqrt{2\lambda}-\epsilon)|x|}\partial^\alpha f \in
L^\infty(\R^3),
\end{equation}
for all $\epsilon >0$ and for any multiindex $\alpha \in \N^3$.
Furthermore
\begin{equation}
\omega(f,X\eta)=0, \quad \forall X \in \mathfrak{g}. \label{e:4.3}
\end{equation}
\end{prop}
\begin{proof}
We first show that a unique solution exists, which follows from
$Q(x)\eta \in (\ker L_+)^\perp$. Indeed, it is sufficient to show this
result for for any $Q_{ij}(x)=x_ix_j$ or $Q_0=1$. By \cite[Theorem 4]{len} we know that $\ker L_+ =
\text{span} \{\partial_1 \eta,\partial_2 \eta,
\partial_3 \eta \}$. Clearly $\la \D_j \eta, \eta \ra = 0$ for all $j \in \{1,2,3\}$. It remains only to show for all $i,j,k \in \{1,2,3\}$ that
\begin{equation} \label{ortho}
\la -\partial_i \eta, x_jx_k \eta\ra =0.
\end{equation}
If $i \neq j$ and $i \ne k$ then \eqref{ortho} is clear, because the integrand is odd in the $x_i$ direction. So we assume $i=j$. If $j \neq k$ then
\[ \la -\partial_i \eta, x_ix_k \eta\ra = -\int \partial_i \eta (x_ix_k) \eta = \int x_k \eta^2 + \int \partial_i \eta (x_ix_k) \eta.\]
But $x_k \eta^2$ is odd in the $x_k$ direction, leading to \eqref{ortho}. A similar argument gives
\eqref{ortho} for $j=k$.

It follows from the PDE solved by $f$ that if $f \in H^s(\R^3)$ then
$f \in H^{s+2}(\R^3)$, implying that $f \in C^\infty(\R^3)$. The
proof of \eqref{expdecayf} now follows closely the proof of
Proposition \ref{prop:etaexpdecay}, and we give it only in outline. We put $w = e^\phi f$ and
introduce
\[
L_+^\phi w \stackrel{\textrm{def}}= e^\phi L_+ e^{-\phi} w = (P_\phi
+ \lambda)w - 2 e^\phi \eta(|x|^{-1}*(\eta e^{-\phi} w)).
\]
We now have
\[
\la L_+^\phi w, w\ra = \frac 12 \int |\nabla w|^2 + \int\left(\wt V - \frac 12|\nabla
\phi|^2 + \lambda\right)w^2 - 2\int e^\phi \eta \left(|x|^{-1}*(\eta
f)\right)w  + \int e^\phi Q(x) \eta w.
\]
Then
\[
\eps \int w^2 \le \int \left(\lambda - \frac 12|\nabla\phi|^2\right)w^2 \le -
\int\wt V w^2 - 2\int e^\phi \eta \left(|x|^{-1}*(\eta f)\right)w  +
\int e^\phi P(x) \eta w.
\]

The $\wt V$ term is handled as before. The two $e^\phi$ factors in
the last term can be absorbed by the $\eta$ factor provided the
exponential growth in $\phi$ is no more than
$e^{\frac{\sqrt{2\lambda}-\eps}2|x|}$. For the middle term, observe that, as in
the case of $\wt V$, the convolution $|x|^{-1}*(\eta f)$ is
continuous and decaying to zero at infinity. Then, the two $e^\phi$
factors can be absorbed by the $\eta$ factor just as in the case of
the last term. In this way we show that
\[ \int w^2 \le C,\]
and proceed as in the proof of Proposition \ref{prop:etaexpdecay}.

We now prove \eqref{e:4.3}. First of all, since $f$ is real,
$\omega(f,e_j \eta)= \im \int f e_j \eta = 0$ for $j \in
\{1,2,3,8\}$ since then $e_j \eta$ is real. Next write
\[f = f_0 + \sum_{j,k=1}^3 f_{jk}, \qquad L_+ f = a_0, \quad L_+f_{jk} = a_{jk}x_jx_k.\]
Since $L_+$ preserves symmetry in $x_k$ for all $k$, we observe that if $j \in \{4,5,6\}$, then
\[
\omega(f_{k\ell},e_j \eta)= \int f_{k\ell} x_{j-1}  \eta = 0,
\]
as the integrand will be odd in some $x_i$ direction.  Finally a calculation
shows that $L_+((2+x \cdot \nabla) \eta) = \eta$, from which it follows that
\[ \omega(f, e_7 \eta) = \int f \eta = \int L_+(f) (2+x \cdot \nabla) \eta = \int (Q(x)\eta)(2+x \cdot \nabla) \eta =0.\]
which completes the proof.
\end{proof}

\section{Reparametrized evolution and proof of Theorem \ref{thm:main2}}\label{reparam}

We write
\[
u(t) = g(t) \cdot(\eta + w(t)), \qquad \omega(w(t),X \eta) = 0 \quad
\forall X \in \mathfrak{g}.
\]
To see that this decomposition is possible, initially for small
times, we apply \ref{3.1}, which allows us to define
\[
g(t) \stackrel{\textrm{def}}= g(u(t)), \quad \tilde u
\stackrel{\textrm{def}}= g(t)^{-1} u(t), \quad w(t)
\stackrel{\textrm{def}}= \tilde u - \eta,
\]
and derive an equation for $w(t)$. Before doing so, however, we introduce some abbreviated notations. For $g(t)$ we write $g = (a,v,\gamma,\mu)$, and 
observe that as a result of $\re \la w, \eta \ra = 0$ and the $L^2$ conservation of the original equation we have
\[2 + \|w\|^2_{L^2} = \|\eta + w\|^2_{L^2} = \|g^{-1} u\|_{L^2}^2 = \mu^{-1} \|u_0\|_{L^2}^2,\]
and hence
\begin{equation}\label{mu}\frac{2 - \eps}{2 + \|w\|^2_{L^2}} \le \mu \le \frac{2 + \eps}{2 + \|w\|^2_{L^2}}, \end{equation}
with $\eps$ as in the statement of Theorem \ref{thm:main2}. This gives a precise sense in which  $\mu \approx 1$. For the remainder of the section we will assume $0 \le \eps \le 1$, although in our theorems $\eps$ is required to be much smaller.

Next we define
\begin{align*}
\alpha &= \alpha(a,\mu) \stackrel{\textrm{def}}= \frac 12 \int V \left(\frac x \mu + a\right)
\eta^2(x)dx - \frac 1{2\mu}\int x\cdot\nabla V\left(\frac x \mu + a\right)\eta^2(x) dx,\\
\beta &= \beta(a,\mu) \stackrel{\textrm{def}}= \frac 1 {2\mu} \int \nabla V\left(\frac x \mu + a\right)\eta^2(x)dx, \\
X &= \mu\sum_{j=1}^3(-\dot a_j + v_j)e_j + \sum_{j=1}^3 \left(\frac {\dot v_j}\mu - \beta_j\right)e_{j+3} +
\left(-\dot\gamma + \dot a \cdot v - \frac 12 |v|^2 + \lambda \mu^2 - \alpha\right)e_7 - \frac {\dot \mu}\mu e_8.
\end{align*}
Observe that $\alpha$ takes values in $\R$, $\beta$ in $\R^3$, and $X$ in $\mathfrak{g}$. Set further
\begin{align*}
\mathcal{L}w &\stackrel{\textrm{def}}= - \frac 12 \Delta w -
\left(|x|^{-1}*\eta^2\right)w - \left(|x|^{-1}*(\eta(w+\bar
w))\right)\eta + \lambda w, \\
\mathcal{N}w &\stackrel{\textrm{def}}= \left(|x|^{-1} *
|w|^2\right)\eta + \left(|x|^{-1} *\eta(w + \bar w)\right)w +
\left(|x|^{-1}*|w|^2\right)w.
\end{align*}
These terms come from writing out $i\Xi_{H}(\eta + w)$. The operator
$\mathcal{L}$ collects the linear terms, and $\mathcal{N}$ the
nonlinear terms.

\begin{lem}\label{weq}
In the above notation, the equation for $w$ is
\begin{align*}
\D_t w &= X \eta + i\left[-V\left(\frac x \mu + a \right) + \alpha + \beta \cdot x\right]\eta \\
&+ X w + i \left[-V\left(\frac x \mu + a\right) + \alpha + \beta
\cdot x\right]w + i\mu^2 \left(-\mathcal{L} + \mathcal{N}\right)w.
\end{align*}
\end{lem}

\begin{proof}The proof of this lemma is a straightforward calculation which follows nearly the same lines as that of \cite[Lemma 3.2]{hz2}, and here we give only a sketch. We first use the definition of $w$ and the chain rule to write
\[\D_t w = - Y(\eta + w) + g^{-1}\Xi_H g(\eta + w),\]
with $Y$ taken from Lemma \ref{yeq}. Next we use Lemma \ref{lem:2.1} to write $g^{-1}\Xi_H g = \mu^{-1} \Xi_{g^*H}$, and compute $\Xi_{g^*H}$ from formula \eqref{e:Hfield2}. Finally, using the soliton equation 
\[-\lambda \eta + \frac 12 \Delta \eta + \left(\frac1 {|x|} * \eta^2 \right)\eta  = 0\]
gives the desired formula.\end{proof}

We now explain the reasons for this notation. Note that if $X=0$, then
\[\dot a = \dot v, \qquad \dot v = - \mu \beta, \qquad  \dot \gamma = \frac 12 |v|^2 + \lambda \mu^2 - \alpha, \qquad \dot \mu=0.\]
giving the equations of motion in \eqref{odes}. In this section and the following section we prove that $|X|$ and $\|w\|_{H^1_x}$ are small, giving Theorem \ref{thm:main2}. Then in Section \ref{final} we give the improvement to Theorem \ref{thm:main} under the necessary additional assumptions on the initial data.

To understand the other crucial features of the notation in Lemma \ref{weq}, we introduce the symplectic projection $P$, characterized by
\[\omega(u,Y\eta) = \omega(P(u)\eta,Y\eta), \quad \forall Y \in \mathfrak{g}.\]
This is given explicitly by
\begin{align*}
P&=\sum_{j=1}^8 e_j P_j, \qquad P_j\colon \mathcal{S}'\to\R \\
P_j(u) &= - \frac 2 {\|\eta\|^2_{L^2}} \omega(u,e_{j+3} \eta) = \re \int u(x) x_j \eta(x) dx, \qquad j \in \{1,2,3\} \\
P_j(u) &= \frac 2 {\|\eta\|^2_{L^2}} \omega(u,e_{j-3}\eta) = - \im \int u(x) \D_{j-3} \eta(x) dx, \qquad j \in \{4,5,6\}\\
P_7(u) &= \frac 2 {\|\eta\|^2_{L^2}} \omega(u,e_8 \eta) = \im \int u(x)(2 + x \cdot \nabla)\eta(x) dx, \\
P_8(u) &= - \frac 2 {\|\eta\|^2_{L^2}} \omega(u, e_7\eta) = \re \int u(x) \eta(x) dx.
\end{align*}

We now compute
\begin{align*}
P(if&(x)\eta(x)) = \sum_{j=4}^6 P_j(if(x)\eta(x)) e_j +
P_7(if(x)\eta(x))e_7 \\
&=-\sum_{j=4}^6 \left(\int
f(x)\eta(x)\D_{j-3}\eta(x)dx\right) e_j + \left(\int f(x)\eta(x)(2 +
x\cdot \nabla)\eta(x)dx\right)e_7\\
&=\frac 1 2\left[-\sum_{j=4}^6 \left(\int
f(x)\D_{j-3}\eta^2(x)dx\right) e_j + \left(\int f(x)\left(4\eta^2(x)
+ x\cdot\nabla \eta^2(x)\right) dx\right)e_7\right]\\
&=\frac 1 2\left[\sum_{j=4}^6 \left(\int
\D_{j-3}f(x)\eta^2(x)dx\right) e_j + \left(\int \left(f(x) -
x\cdot\nabla f(x)\right)\eta^2(x) dx\right)e_7\right]\\
&\stackrel{\textrm{def}}= i \alpha + i\beta \cdot x.\\
\end{align*}
Observe that in the case that $f(x) = V(x/\mu + a)$ these $\alpha$
and $\beta$ agree with those defined previously.

We have the following Taylor expansions, where $\delta_{jk}$ is the Kronecker delta:
\begin{align*}
V\left(\frac x \mu + a \right) &= V(a) + \nabla V(a) \cdot \frac x
\mu + \frac 1 {\mu^2} \sum_{j,k=1}^3 \left(1 - \frac {\delta_{jk}}2\right)x_jx_k \D_j\D_k V(a)  + \Oh(h^3),\\
\alpha &= V(a) + \frac 3 {4\mu^2} \int \left[\sum_{j=1}^3x_j^2\D_j^2 V(a)\right] \eta^2(x)dx + \Oh(h^3),\\
\beta&= \frac{\nabla V(a)}\mu + \Oh(h^3),
\end{align*}
and thus
\begin{align*}
-V&\left(\frac x \mu + a \right) + \alpha + \beta \cdot x \\
&=- \frac 1 {\mu^2} \sum_{j,k=1}^3 \left(1 - \frac {\delta_{jk}}2\right)x_jx_k \D_j\D_k V(a) + \frac 3 {4\mu^2} \int \left[\sum_{j=1}^3x_j^2\D_j^2 V(a)\right] \eta^2(x)dx+ \Oh(h^3),\\
&\stackrel{\textrm{def}}= \sum_{j,k=1}^3a_{jk}x_jx_k + a_0+ \Oh(h^3) \stackrel{\textrm{def}}= Q(x) +  \Oh(h^3).
\end{align*}
where all the errors are polynomially bounded in $x$. In the sequel we will apply Proposition \eqref{4.2} using this $Q(x)$. Observe that it satisfies the necessary orthogonality condition because $\omega(i(V(x/ \mu + a),X\eta)) = 0$, and $Q(x)$ is of order $h^2$.

We now study $w$ by writing $w = \tilde w + w_1$, where $\tilde w$ solves away the principal forcing terms of the equation of $w$. More precisely, we put 
\begin{align*}
\tilde w &\stackrel{\textrm{def}}= \sum_{j,k = 1}^3 \tilde w_{jk},
\qquad \tilde w_{jk} \stackrel{\textrm{def}}= - \frac {\D_j\D_k
V(a)} {\mu^4} f_{jk}, \\ f_{jk} &\stackrel{\textrm{def}}=
L_+^{-1}\left(- \sum_{j,k=1}^3 \left(1 - \frac 12 \delta_{jk}\right)x_jx_k  + \delta_{jk}\frac 3 4\int x_j^2\eta^2(x)dx\right)\eta.
\end{align*}

Then $\tilde w$ satisfies the PDE
\begin{multline*}
\D_t \tilde w = - i \mu^2 \mathcal{L}\tilde w - \frac i{\mu^2}\left(- \sum_{j,k=1}^3 \left(1 - \frac 12 \delta_{jk}\right) x_jx_k\D_j\D_kV(a)  + \frac 34
\int \left[\sum_{j=1}^3 x_j^2\D_j^2V(a)\right]\eta^2(x)dx\right)\eta \\
+ \sum_{j,k=1}^3\theta_{jk} f_{jk},
\end{multline*}
where
\[
\theta_{jk}(t) \stackrel{\textrm{def}}= \frac{d}{dt}\left[\frac{-\D_j\D_k V(a)}{\mu^4}\right] = \frac{-\D_j\D_k \nabla V(a)\cdot \dot
a}{\mu^4} + \frac{4\D_j\D_k V(a) \dot \mu}{\mu^5}
\]

\begin{lem}
There exists an absolute constant $c$ such that if $\|w\|_{H^1} \le 1/c$, then
\[|X| \le c (h^2\|w\|_{H^1} + \|w\|^2_{H^1} + \|w\|^3_{H^1}).\]
\end{lem}

\begin{proof}
Since $P w_t = \D_t Pw = 0$, Lemma \ref{weq} gives
\begin{multline*}X = P(i(V(x/\mu + a) - \alpha - \beta \cdot x)\eta) + P(i(V(x/\mu + a) - \alpha - \beta \cdot x)w) - P(Xw)\\ - \mu^2P(i\mathcal{N}w) - \mu^2P(i\mathcal{L}w).\end{multline*}
We have already observed that the first term vanishes. Next the estimate $|P(Yw)|\le c |Y|\|w\|_{H^1}$ shows that
\[|P(i(V(x/\mu + a) - \alpha - \beta \cdot x)w)| \le c h^2\|w\|_{H^1}, \qquad |P(Xw)| \le c|X|\|w\|_{H^1}.\]

For the $P(i\mathcal{N}w)$ term we must estimate the following integral, where
$\psi_k$ are taken from $w,\eta, e_j\eta$,:
\begin{align}
\label{convo}\int\left|\left(|x|^{-1}*(\psi_1\psi_2)\right)\psi_3
\psi_4\right| &\le
\||x|^{-1}*(\psi_1\psi_2)\|_{L^3}\|\psi_3\|_{L^6}\|\psi_4\|_{L^2}
\notag\\
&\le c\|\psi_1\psi_2\|_{L^1}\|\psi_3\|_{L^6}\|\psi_4\|_{L^2} \le
c\|\psi_1\|_{L^2}\|\psi_2\|_{L^2}\|\psi_3\|_{H^1}\|\psi_4\|_{L^2}
\end{align}
For this we have used H\"older's inequality, the Hardy-Littlewood-Sobolev inequality, and Sobolev embedding. This results in
\[|P(i\mathcal{N}w)| \le c(\|w\|_{H^1}^2 + \|w\|_{H^1}^3).\]

Finally, from \eqref{il} we have
\[ P(i\mathcal{L}w) = 0,\]
which combines with the previous estimates to give
\[|X| \le c h^2 \|w\|_{H^1} + c|X|\|w\|_{H^1} + c(\|w\|_{H^1}^2 + \|w\|_{H^1}^3).\]
Here we have removed the factors of $\mu$ using \eqref{mu}. If $\|w\|_{H^1}$ is sufficiently small, this implies the desired inequality.
\end{proof}

\begin{lem}\label{5.1}
Suppose there are positive constants $c_1$, and $h_0$ such that
\[
\|w\|_{L^\infty_{[t_1,t_2]}H^1_x} \le c_1 h^{\frac 12 + \delta}, \qquad h^{2+2\delta}(t_2 - t_1)\la t_2 - t_1 \ra \le c_1, \qquad 0 < h \le h_0,
\]
for some $t_1 < t_2$, $\delta \ge 0$.
Then
\[
\sup_{t_1 < t < t_2} |\theta(t)| \le ch^3, \qquad \sup_{t_1 < t < t_2} |v(t)| \le c,\]
for a constant $c$ depending only on $c_1$, $h_0$,
$\|W\|_{C^3(\R^3)}$ and $|v(t_1)|$.
\end{lem}

\begin{proof}
The conclusion concerning $\theta$ will follow from $|\dot\mu| \le ch^{1 + 2\delta}$ and $|\dot a| \le c$. Observe that our assumption on $w$ implies that the bounds for $\mu$ in \eqref{mu} can be improved to
\[1 - ch^{\frac 12 + \delta} \le \mu \le 1 + ch^{\frac 12 + \delta}.\]
By the definition of $X$ and the Taylor expansions and the bound on $X$, we have
\[
\left|\frac{\dot v}\mu + \nabla V(a)\right| + \left|\frac{\dot \mu}\mu\right| + |\mu(-\dot a + v)| \le c |X| \le c(h^2\|w\|_{H^1} + \|w\|_{H^1}^2 + \|w\|_{H^1}^3),
\]
which immediately gives the desired bound on $|\dot\mu|$. For the
bound on $|\dot a|$, it suffices to prove $|v| \le c$, which we do
by first integrating the above inequality to obtain:
\[
\label{5.7} \sup_{t_1 < t < t_2} |v(t)| \le |v(t_1)| + ch\|\nabla
W\|_{L^\infty}(t_2-t_1) + c|X|(t_2-t_1).
\]
Next we prove a near conservation of classical energy:
\begin{align*}
\sup_{t_1\le t \le t_2} \Big|\left(\frac {|v|^2} 2 + V(a)\right)& -
\left(\frac {|v(t_1)|^2} 2 + V(a(t_1))\right)\Big| \\
&\le (t_2 - t_1)\sup_{t_1\le t \le t_2} |\dot v \cdot v + \nabla V
\cdot a| \\
&\le (t_2 - t_1)\sup_{t_1\le t \le t_2} \left( |\dot v + \nabla
V(a)||v| + |\nabla V(a)||\dot a - v| \right)\\
&\le c(t_2 - t_1)\left[|X|\sup_{t_1\le t \le t_2}|v| +
h\|\nabla W\|_{L^\infty}|X|\right]\\ &\le c|X|(t_2 - t_1)\left[|v(t_1)| + ch\|\nabla W\|_{L^\infty}\la
t_2-t_1\ra + c|X|(t_2-t_1)\right].
\end{align*}
From this it follows that $\sup_{t_1 \le t \le t_2} |v(t)| \le c$, which concludes the proof.
\end{proof}

This will be crucial for the estimate of the true error $w$.

\begin{lem}[Lyapounov energy estimate]\label{5.2}
Suppose that, for some constants $c_1$ and $h_0$,
\[\|w\|_{L^\infty_{[t_1,t_2]}H^1_x} \le c_1 h^{\frac 12}, \qquad 0< h \le h_0. \]
Then, provided
\[|t_2 - t_1| \le \frac {c_2} h,\]
we have
\[
\|w\|_{L^\infty_{[t_1,t_2]}H^1_x} \le  c_3\|w_1(t_1)\|_{H^1} + c_4h^2.
\]
The constants $c_2$ and $c_4$ depend only upon $c_1$, $h_0$,
$\|W\|_{C^3(\R^3)}$ and $|v(t_1)|$.
The constant $c_3$ is an absolute constant.
\end{lem}

We postpone the proof of this lemma to the end of the section, first
demonstrating how it is applied in the bootstrap argument. We prove the
following proposition, from which Theorem \ref{thm:main2} follows.

\begin{prop}
Let $w_0 = w(0)$ and fix constants $\tilde c_1 >0$ and $\delta_0 \in (0,3/4)$.
Then there exist constants $h_0$ and $c$ such that if
\[ 0 \le \delta \le \delta_0, \qquad 0 < h \le h_0, \qquad\|w_0\|_{H^1} \le \tilde c_1 h^{\frac 12 + 3 \delta_0}, \qquad0 < T \le \frac{\tilde c_1}h + \frac {\delta|\log h|}{ch}
\]
then
\[
\|w\|_{L^\infty_{[0,T]}H^1_x} \le c h^{-\delta}\left(\|w_0\|_{H^1} + h^2\right).
\]
The constants $h_0$ and $c$ depend only on $\tilde c_1$, $\delta_0$, $|v(0)|$, and $\|W\|_{C^3(\R^3)}$.
\end{prop}

\begin{proof}
To apply Lemma \ref{5.2}, we observe that by the continuity in $t$ of $\|w\|_{L^\infty_{[0,t]}H^1_x}$ we know immediately that the hypotheses are satisfied on $[0,t]$ for sufficiently small $t$. At this point the conclusion of the lemma tells us that at the end of this interval the error is still small enough that we may proceed for larger $t$, until we reach $t = c_2/h$. In this way we apply Lemma \ref{5.2}, $k$ times on successive intervals of
length $c_2/h$, where $c_2$ and $k$ will be fixed later, giving the bound
\[
\|w\|_{L^\infty_{[0,c_2k/h]}H^1_x} \le c_3^k \|w_0\|_{H^1} + \left(\sum_{j=0}^{k-1} c_3^j\right) c_4 h^2.
\]
This is only valid provided that the hypotheses of Lemmas \ref{5.1} and \ref{5.2} are satisfied over the whole collection of time intervals. We must use Lemma \ref{5.1} to control $|v|$ uniformly over the full time interval $[0,c_2 k/h]$, and to apply this we need
\[c_3^k \|w_0\|_{H^1} + \left(\sum_{j=0}^{k-1} c_3^j\right) c_4 h^2 \le c_1 h^{\frac 12 + \delta}, \qquad c_2^2k^2 h^{2\delta} \le c_1,\]
for some constant $c_1$. We will determine $c_1$ momentarily, and at that point $c_2$ will be the constant which emerges from Lemma \ref{5.2}. If
\[k = \frac{\tilde c_1}{c_2} + \delta \frac{|\log h|}{\log c_3},\]
it suffices to have
\begin{equation}\label{hypo}c_3^{\tilde c_1/c_2}\tilde c_1 h^{\frac 12 + 3\delta_0 - \delta} + c_3^{\tilde c_1/c_2} c_4 h^{2-\delta}\le c_1 h^{\frac 12 + \delta}, \qquad \tilde c_1^2 \left\la\delta \frac{|\log h|}{\log c_3}\right\ra^2 h^{2\delta} \le c_1.\end{equation}

We are now ready to choose our constants. We first take $c_1$ such that the second inequality of \eqref{hypo} holds. Then $c_2$ is given by Lemma \ref{5.2}, and we take $h_0$ is such that the first inequality of \eqref{hypo} holds. Note that the hypotheses of Lemma \ref{5.1} are satisfied a fortiori.
\end{proof}

It now remains only to prove Lemma \ref{5.2}.

\begin{proof}[Proof of Lemma \ref{5.2}]
In this proof, unless otherwise mentioned, all constants depend only
upon $c_1$, $\|W\|_{W^{\infty,3}}$ and $|v(t_1)|$.

Let
\[ w_1 \stackrel{\textrm{def}}= w - \tilde w,\]
Now
\begin{align*}
\D_t &w_1 = - i\mu^2\mathcal{L} w_1 + X \eta - \theta f\\
&+i \left[-V\left(\frac x \mu + a \right) + \alpha + \beta \cdot x -
\frac x {2\mu^2} \cdot\nabla^2 V(a) x  + \frac 3
{2\mu^2\|\eta\|^2_{L^2}} \int x \cdot \nabla^2V(a)x
\eta^2(x)dx\right] \eta \\
&+Xw +i\left[-V\left(\frac x \mu + a \right) + \alpha + \beta \cdot
x\right]w + i\mu^2\mathcal{N} w.
\end{align*}
By grouping forcing terms into $f_1$, we rewrite the above as
\[
\D_t w_1 = - i\mu^2\mathcal{L} w_1 + X \eta + f_1 +Xw +
i\left[-V\left(\frac x \mu + a \right) + \alpha + \beta \cdot
x\right]w + i\mu^2\mathcal{N} w,
\]
observing that, using Lemma \ref{5.1}, we have $\|f_1\|_{H^1} \le ch^3$

We recall that $\mathcal{L}$ is self-adjoint with respect to
\[\la u, v \ra = \re\int u \bar v,\]
and hence
\begin{align*}
\frac 12 \D_t \la \mathcal{L}&w_1,w_1\ra = \la \mathcal{L} w_1, \D_t
w_1 \ra \\
= &- \mu^2\la\mathcal{L} w_1, i\mathcal{L} w_1\ra + \la \mathcal{L} w_1, X\eta \ra +
\la\mathcal{L} w_1,f_1\ra + \la\mathcal{L} w_1,Xw_1\ra
+\la\mathcal{L} w_1,X\tilde w\ra \\
&+ \la\mathcal{L} w_1,i\left[-V\left(\frac x \mu + a \right) +
\alpha + \beta \cdot x\right]w_1\ra + \la\mathcal{L}
w_1,i\left[-V\left(\frac x \mu + a \right) + \alpha + \beta \cdot
x\right]\tilde w\ra \\
&+  \la\mathcal{L} w_1,i\mu^2\mathcal{N} w\ra \\
 = & \,\,\,\textrm{I
+ II + III + IV + V + VI + VII + VIII}
\end{align*}

Now we analyze these terms one-by-one. First
\[
\textrm{I} = \textrm{II} = 0.
\]
In the case of I this follows from \eqref{e:products}, the definition of $\la\cdot,\cdot\ra$. In the case of II, we recall that $\omega(w,X\eta) = 0$ by construction of $w$, and that $\omega(\tilde w, X \eta) = 0$ from \ref{e:4.3}, as a result of which we have $\omega(w_1, X \eta) = 0$. Finally $\omega(i\mathcal{L} w_1,X\eta) = 0$ by \eqref{il}, and then we use \eqref{e:products} to relate $\la\cdot,\cdot\ra$ and $\omega(\cdot,\cdot)$.

Next we show that
\[|\textrm{III}| \le c \|w_1\|_{H^1} \|f_1\|_{H^1} \le ch^3\|w_1\|_{H^1}.\]
This estimate is straightforward in the case of the convolution-free
terms of $\mathcal{L}$. For the terms with convolutions, we apply
\eqref{convo} with $f_1$ in place of $\psi_4$ and the other $\psi_k$
chosen appropriately from among $\eta, w, \bar w$.

Next we look at $\textrm{IV} = \la \mathcal{L} w_1, X w_1\ra$. We
first recall that $X = \sum_{j=1}^8 a_j e_j$ with $|a_j| \le
c(h^2\|w\| + \|w\|_{H^1}^2 + \|w\|_{H^1}^3)$. We the proceed term by term
according to $\mathcal{L}w_1 = \frac 12 w_1 - \frac 12 \Delta w_1 -
\left(|x|^{-1}*\eta^2\right)w_1 - \eta\left(|x|^{-1}*(\eta(w_1+\bar
w_1))\right)$:
\begin{align*}
\la w_1, X w_1 \ra &= a_8\la w_1, 2w_1 + x \cdot \nabla w_1\ra = 
\frac 12 a_8 \la w_1, w_1 \ra,\\
\la \Delta w_1, X w_1\ra &= \sum_{j=1}^3 a_{j+3}\la \Delta w_1, i
x_j w_1\ra + a_8\la \Delta w_1, 2w_1 + x \cdot \nabla w_1\ra \\
&= \sum_{j=1}^3 a_{j+3}\la \D_j w_1,  i w_1\ra + \frac 12 a_8\la
\nabla w_1,\nabla w_1\ra,
\end{align*}
and thus the above two terms are bounded by $c |X|\|w_1\|_{H^1}^2$. For
the terms involving $\eta$ we use \eqref{convo} to obtain the same
bound, giving
\[|\textrm{IV}| \le c(h^2 + \|w\|_{H^1} + \|w\|_{H^1}^2)\|w_1\|_{H^1}^3.\]

Next $\textrm{V} = \la\mathcal{L} w_1,X\tilde w\ra$ has a similar
expansion, but including more nonzero terms. We estimate these terms
as before in \eqref{convo}, using H\"older's inequalty,
Hardy-Littlewood-Sobolev, and Sobolev embedding, to obtain
\[
|\textrm{V}| \le c|X| \|w_1\|_{H^1}\|\la x \ra
\tilde w \|_{H^2}.
\]
However, $\|\la x \ra \tilde w \|_{H^2} \le c h^2$, giving
\[ |\textrm{V}| \le ch^2(h^2 + \|w\|_{H^1} + \|w\|_{H^1}^2)\|w_1\|_{H^1}. \]

For $\textrm{VI}$ once again we obtain a number of vanishing terms:
\begin{align*}
\textrm{VI} &= \la \mathcal{L} w_1, i \left[-V\left(\frac x \mu + a
\right) + \alpha + \beta \cdot x \right]w_1\ra \\
&= \la - \frac 12 \Delta w_1 - \eta\left(|x|^{-1} * (\eta(w_1 + \bar
w_1))\right), i \left[-V\left(\frac x \mu + a \right) + \alpha +
\beta \cdot x \right]w_1\ra.
\end{align*}
To estimate the first term, we integrate by parts as before and use
\[\left| - \frac 1 \mu \nabla V\left(\frac x \mu + a \right) + \beta \right| \le c h.\]
For the second term, we use \eqref{convo} together with
\[
\left|\left[-V\left(\frac x \mu + a \right) + \alpha + \beta \cdot x
\right]\eta\right|\le c h^2.
\]
This gives the bound
\[ |\textrm{VI}| \le c h \|w_1\|_{H^1}^2. \]

For $\textrm{VII}$ we proceed in the same way, without the vanishing
terms but also without the restriction that only $H^1$ norms may be
used. We obtain
\begin{align*}
|\textrm{VII}| &\le c\|w_1\|_{H^1}\left\|\left[-V\left(\frac x \mu +
a \right) + \alpha + \beta \cdot x \right]\tilde w\right\|_{H^1} \\
&\le c h^2\|w_1\|_{H^1}\|\la x \ra^2 \tilde w\|_{H^1} \le c
h^4\|w_1\|_{H^1}.
\end{align*}

Finally, for $\textrm{VIII} = \la\mathcal{L} w_1,i\mu^2\mathcal{N}
w\ra$ we write $w = w_1 + \tilde w$ and expand. We integrate by
parts for the $\Delta$ term, and use \eqref{convo}, twice as needed
for the terms with two convolutions. This allows us to put all
factors in an $H^1$ norm, giving a bound of
\[
|\textrm{VIII}| \le c\left(h^6 \|w_1\|_{H^1} + h^4 \|w_1\|^2_{H^1} +
h^2 \|w_1\|^3_{H^1} + \|w_1\|^4_{H^1} \right)
\]

Combining all this gives
\[
|\D_t \la \mathcal{L} w_1, w_1 \ra | \le c\left(h^3 \|w\|_{H^1} +
h \|w\|_{H^1}^2 + h^2 \|w\|_{H^1}^3 + \|w\|_{H^1}^4 + \|w\|_{H^1}^5 \right).
\]
From Appendix \ref{app:b} we have uniform boundedness of $\|u\|_{H^1}$, while from Lemma \ref{5.1} we have uniform boundedness of $|v|$ over our time interval, from which we conclude that $\|w\|_{H^1} \le c$, and hence
\[
|\D_t \la \mathcal{L} w_1, w_1 \ra | \le c\left(h^3 \|w\|_{H^1} +
h \|w\|_{H^1}^2 + \|w\|_{H^1}^4\right).
\]
Now we use $w = w_1 + \tilde w$ to write $\|w\|_{H^1} \le c(\|w_1\|_{H^1} + h^2)$ and hence
\[
|\D_t \la \mathcal{L} w_1, w_1 \ra | \le c\left(h^5  +
h \|w_1\|_{H^1}^2 +  \|w_1\|_{H^1}^4 \right).
\]
Integrating in time gives
\[
\la \mathcal{L} w_1(t), w_1(t) \ra \le \la \mathcal{L} w_1(t_1),
w_1(t_1)\ra + c(t-t_1)\left(h^5 + h \|w_1\|_{H^1}^2 +
\|w_1\|_{H^1}^4\right)
\]
From \eqref{e:4.1} we have
\[
\|w_1(t)\|^2_{H^1} \le c\la \mathcal{L}w_1(t),w_1(t)\ra,
\]
and by direct esimation we have
\[
|\la \mathcal{L}w_1(t),w_1(t)\ra| \le c \|w_1(t)\|^2_{H^1}.
\]
This leads to
\begin{align*}
\|w_1&\|^2_{L^\infty_{[t_1,t]}H^1_x} \le \tilde
c\|w_1(t_1)\|^2_{H^1} \\
&+ c(t-t_1)\left(h^5 + h
\|w_1\|_{L^\infty_{[t_1,t]}H^1_x}^2 +
\|w_1\|_{L^\infty_{[t_1,t]}H^1_x}^4\right),
\end{align*}
with $\tilde c$ an absolute constant. Requiring that $t_2 - t_1 \le
c_2/h$ for a small constant $c_2$, and subtracting the quadratic
term to the left hand side implies
\[
\|w_1\|^2_{L^\infty_{[t_1,t]}H^1_x} \le 2\tilde
c\|w_1(t_1)\|^2_{H^1} + c(t_2-t_1)\left(h^5 +
h\|w_1\|_{L^\infty_{[t_1,t]}H^1_x}^4\right).
\]
This is a quadratic inequality in
$\|w_1\|^2_{L^\infty_{[t_1,t]}H^1_x}$. In general,
\[A > 0,\, B > 0,\, X \in \R,\, BX^2 - X + A \ge 0,\, X \le (2B)^{-1},\, 4AB < 1 \,\Longrightarrow\, X \le 2A.\]
In our case, assuming that
\[(t_2-t_1)h \|w_1\|^2_{L^\infty_{[t_1,t]}H^1_x} + (t_2-t_1)^2h^6\le c_2\]
we have
\[
\|w_1\|^2_{L^\infty_{[t_1,t_2]}H^1_x} \le 4\tilde
c\|w_1(t_1)\|^2_{H^1} + ch^5(t_2-t_1).
\]
From this, together with $w = w_1 + \tilde w$ the desired result follows.

\end{proof}

\section{Proof of Theorem \ref{thm:main}}\label{final}

In this section we make use of the following lemma:
\begin{lem}\label{6.1}
Suppose that $0 < h \ll 1$, and $a=a(t),v=v(t),
\epsilon_1=\epsilon_1(t), \epsilon_2=\epsilon_2(t)$ are $C^1$ real
valued functions. Suppose $f: \R^3 \rightarrow \R$ is $C^2$ mapping
such that $|f|$ and $|f^\prime|$ are uniformly bounded. Suppose that
on $[0,T]$,
\[
\left\{
\begin{array}{ll}
\dot a = v + \epsilon_1, &a(0)=a_0\\
\dot v = hf(ha)+ \epsilon_2, & v(0)=v_0
\end{array}
\right.
\]
Let $\ol a = \ol a(t)$ and $\ol v = \ol v(t)$ be the $C^1$ real
valued functions satisfying the exact equations
\[
\left\{
\begin{array}{ll}
\dot {\ol a} = \ol v + \epsilon_1, &\ol a(0)=a_0\\
\dot {\ol v} = hf(h \ol a)+ \epsilon_2, & \ol v(0)=v_0
\end{array}
\right.
\]
with the same initial data. Suppose that on $[0,T]$, we have
$|\epsilon_j| \le h^{4-\delta}$ for $j=1,2$. Then provided $T \le
c h^{-1} + \delta h^{-1} \log (1/h)$, we have on $[0,T]$ the estimates
\[
|a-\ol a| \le \tilde c h^{2-2 \delta} \log(1/h), \quad |v - \ol v| \le
\tilde c h^{3-2\delta} \log (1/h)
\]
\end{lem}
The statement and proof of this lemma is almost identical to those of
\cite[Lemma 6.1]{hz2}. The only change in this proof is that we use
$g = \int_0^1 \nabla f(h \ol a +t (ha-h\ol a)) dt$.

For Theorem \ref{thm:main} we assume $\eps = \Oh(h^2)$, in which case we have the following ODEs for $a$ and $v$:
\[\dot a = v + \Oh(h^{4 - 4\delta}), \qquad \dot v = -\frac 12 \int \nabla V(x+a)\eta^2(x)dx +  \Oh(h^{4 - 4\delta}).\]
Lemma \ref{6.1} allows us to replace these with
\[\dot a = v, \qquad \dot v = -\frac 12 \int \nabla V(x+a)\eta^2(x)dx.\]
Direct integration of the error terms in the equations for $\mu$ and $\gamma$ allows them to be dropped as well, giving Theorem \ref{thm:main}.

\appendix

\section{Properties of $\eta$}\label{app:eta}

In this appendix we review the properties of the function $\eta$ which we need in this paper. This material is essentially well-known, and further information and references may be found in \cite{len}. First we recall a lemma
from \cite[Appendix A]{len}.
\begin{lem}\label{lem:len}
For each $\lambda > 0$, the equation
\begin{equation} \label{e:eta}
-\frac{1}{2} \Delta \eta + \wt V \eta = - \lambda\eta
\end{equation}
with $\wt V=-|x|^{-1}*\eta^2$, has a unique radial, nonnegative
solution $\eta \in H^1(\R^3)$ with $\eta \not\equiv 0$. Moreover, we
have that $\eta(r)$ is strictly positive.
\end{lem}

In this paper we choose $\lambda$ such that
\[\|\eta\|_{L^2}^2 = 2.\]

We will also need the following exponential decay result.

\begin{prop}\label{prop:etaexpdecay}
Let $\eta \in H^1(\R^3;\R)$ satisfy \eqref{e:eta}. Then $\eta \in
C^\infty(\R^3)$, and for any multiindex $\alpha$ and $\epsilon >0$
there exists $C$ such that
\[
\left|\partial^\alpha \eta(x)\right| \le C e^{-(\sqrt{2\lambda}-\epsilon)|x|}.
\]
\end{prop}
\begin{proof}
Observe first that $\wt V$ is continuous and obeys $\lim_{|x|
\rightarrow \infty} \wt V =0$. Indeed, write $|x|^{-1}=\chi_1+
\chi_2$, where $\chi_1$ is smooth and agrees with $|x|^{-1}$ near
infinity, and $\chi_2$ is compactly supported and in $L^p$ for
$p<3$. The $\chi_1$ terms is clearly smooth, and we prove the decay
by treating it in two pieces:
\[
\int_{|y|\le |x|/2} \chi_1(x-y)\eta^2(y)dy \le \int_{|y| \le |x|/2}
\frac{C}{\la x-y \ra}\eta^2(y) dy \le \frac{C}{|x|} \| \eta
\|^2_{L^2}
\]
\[
\int_{|y| \ge |x|/2} \chi_1(x-y) \eta^2(y) dy \le
\|\chi_1\|_{L^\infty} \int_{|y| \ge |x|/2} \eta^2(y) dy
\]
On the other hand note that since $\eta \in H^1(\R^3)$, the
Gagliardo-Nirenberg inequality implies that $\eta \in L^6(\R^3)$,
and in particular $\eta^2 \in L^2$. Thus $\chi_2 * \eta^2$ has a Fourier
transform in $L^1$, giving the desired regularity and decay.

Now it follows from \eqref{e:eta} that $\eta \in H^2$.
Differentiating the equation and applying the previous argument shows
that $\eta \in H^3$. By induction we find that $\eta \in H^s$, and in particular $\eta \in C^\infty$.

We now prove the exponential decay as follows. Let
$P=-\frac{1}{2}\Delta + \wt V$, let $\phi \in C^\infty$ be bounded
together with its first derivatives, and let
\[
P_\phi \stackrel{\textrm{def}}= e^\phi P e^{-\phi} = -
\frac{1}{2}\Delta + \nabla \phi \cdot \nabla - \frac{1}{2}|\nabla
\phi|^2 +\frac{1}{2}\Delta \phi+ \wt V.
\]
Let $w = e^\phi \eta$ and , observing that integrating by parts
gives $\int(\nabla \phi \cdot \nabla w)w = - \int (\nabla \phi \cdot
\nabla w)w - \int (\Delta \phi)w^2,$ write
\[
0=\left\la \left(P_\phi +\lambda\right)w,w \right\ra_{L^2} = \frac 12\int | \nabla w |^2 + \int \left(\wt V
+\lambda - \frac 12|\nabla \phi|^2\right)w^2
\]
Now, provided $| \nabla \phi |^2 \le 2\lambda- 2\epsilon$ we have
\[
\epsilon \int w^2 \le \int \left(\lambda-\frac{|\nabla \phi|^2}2\right)w^2 \le - \int \wt V
w^2  \le \frac{\epsilon}{2} \int_{\{x:\wt V(x) \ge
-\epsilon/2\}}w^2 - \int_{\{x:\wt V(x) < -\epsilon/2\}} \wt V w^2.
\]
The integral over $\{x:\wt V(x) \ge - \epsilon/2\}$ can now be
subtracted to the other side of the inequality, while $\{x:\wt V(x)
< - \epsilon/2\}$ is a bounded set as a result of $\lim_{|x|
\rightarrow \infty} \wt V (x)=0$. We may then write
\[
\int w^2 \le C
\]
where $C$ depends on $\eta$, $\sup|\phi|$, and $\epsilon$. If we
apply this result with a sequence of functions $\phi_n$ such that
$\phi_n=(\sqrt{2\lambda - 2\epsilon})x_1$ on the ball of radius $n$ and is modified
outside that ball to be smooth with bounded derivatives, we find
that $e^{\sqrt{2\lambda - 2\epsilon}x_1}\eta \in L^2$, and similarly
\[
e^{\sqrt{2\lambda - 2\epsilon}|x|}\eta(x) \in L^2
\]
Differentiating \eqref{e:eta} and applying the same argument proves
that
\[
e^{\sqrt{2\lambda - 2\epsilon}|x|} \partial^\alpha \eta(x) \in L^2,
\]
from which the desired result follows.
\end{proof}

\section{Well-posedness}\label{app:b}

In this appendix we prove well-posedness for the equation \eqref{e:hartree} in $H^1(\R^3)$. This result is known (see for example \cite{c}), but for the reader's convenience we review the result in the special case which we study here. We adopt the notation $\|u\|_{W^{k,p}} = \sum_{|\alpha|\le k} \|\D^\alpha u\|_{L^p}$.

We will use the following Strichartz estimates (see for example \cite{kt}). 
\begin{lem} Suppose $q,r,\tilde q', \tilde r' \in [1,\infty]$ satisfy
\[\frac 2 q + \frac n r = \frac n 2, \qquad \frac 2 {\tilde q'} + \frac n {\tilde r'} = \frac {4+n} 2.\]
Then 
\[\|e^{it\Delta}u_0\|_{L^q_{[0,T]}L^r_x} \le c \|u_0\|_{L^2} \qquad \left\|\int_0^t e^{i(t-s)\Delta} f(s) ds\right\|_{L^q_{[0,T]}L^r_x} \le c\|f\| _{L^{\tilde q'}_{[0,T]}L^{\tilde r'}_x},\]
for all $u_0 \in L^2(\R^n)$ and $f \in L^{\tilde q'}([0,T],L^{\tilde r'}(\R^n))$.
\end{lem}

In the remainder of this section only, $c$ denotes a constant which may vary from line to line, but is absolute, that is independent of all parameters in the problem. Let $V \in W^{1,\infty}(\R^3,\R)$, and let $u_0 \in H^1(\R^3)$ be
given, and define
\[
N(u) = -\left(|x|^{-1}*|u|^2\right)u, \qquad F(u)(t) =
e^{it\Delta}u_0 - i \int_0^t e^{i(t-s)\Delta} \left[N(u(s)) + Vu(s)\right]ds.
\]

A function $u$ solves the Hartree equation if and only if it is a fixed point of $F$. We have the following

\begin{lem} For any $T>0$, we have\
\begin{align*}
\|N(u)\|_{H^1(\R^3)} &\le c\|u\|_{L^2(\R^3)} \|\nabla
u\|_{H^1(\R^3)}, \\
\|F(u)\|_{L^\infty([0,T],H^1(\R^3))} &\le \|u_0  \|_{H^1(\R^3)} +
T^{1/2}(c\|u\|^3_{H^1(\R^3)} +
\|V\|_{W^{1,\infty}(\R^3)}\|u\|_{H^1(\R^3)}),
\end{align*}
where $c$ is an absolute constant.
\end{lem}

\begin{proof}
We first compute
\begin{equation}\label{appb}
\left\|\left(|x|^{-1}*|u|^2\right)u\right\|_{L^2} \le
\left\|\left(|x|^{-1}*|u|^2\right)\right\|_{L^3}\|u\|_{L^6} \le
c\left\||u|^2\right\|_{L^{1}}\|u\|_{L^6} \le c\|\nabla
u\|_{L^2}\|u\|^2_{L^2},
\end{equation}
where we have used in the first inequality H\"older, in the second
Hardy-Littlewood-Sobolev, and in the third H\"older followed by the
Sobolev inclusion $\dot{H}^1(\R^3) \subset L^6(\R^3)$. From this the result concerning $N$ follows.

We now look at $F$. We have
$\|e^{it\Delta}u_0\|_{L^\infty([0,T],H^1(\R^3))} =
\|u_0\|_{H^1(\R^3)}$ because the Schr\"odinger propagator is unitary
on all Sobolev spaces. We then compute using Strichartz estimates
that
\[
\left\|\int_0^t e^{i(t-s)\Delta} N(u(s))
ds\right\|_{L^\infty([0,T],L^2(\R^3))} \le c\|N(u)\|_{L^2_{[0,T]} L^{6/5}_x} \le cT^{1/2} \|N(u)\|_{L^\infty_{[0,T]} L^{6/5}_x}
\]
% T \|N(u(t))\|_{L^\infty([0,T],H^1(\R^3))},
Using the same sequence of inequalities as in \eqref{appb} we get that
\[\left\| \left( |x|^{-1} * |u|^2\right)u \right\|_{L^{6/5}} \le \left\| |x|^{-1} * |u|^2\right\|_{L^3} \left\| u
\right\|_{L^2} \le c\| |u|^2\|_{L^1} \|u\|_{L^2} = c \|u\|_{L^2}^3\]
The same arguments show that
\[\left\|\nabla \int_0^t e^{i(t-s)\Delta} N(u(s)) ds \right\|_{L^\infty([0,T],L^2(\R^3))} \le T^{1/2} \|u\|_{L^2}^2 \| \nabla u\|_{L^2}.\]
The result concerning $F$ follows from this.
\end{proof}

%Using this we prove local existence of a solution to Hartree
%equation
%\begin{equation} \label{e:htnov}
%i \D_t u = - \frac 12 \Delta u - \left(|u|^2 * \frac 1 {|x|}\right)u
%\\ u(x,0) = u_0(x) \in H^1(\R^3;\C).
%\end{equation}
\begin{prop}\label{p:locwp}
For each $u_0 \in H^1(\R^3;\C)$ there exists $T \in \R$ such that
\eqref{e:hartree}
 has a solution $u(x,t) \in L^\infty([0,T],H^1(\R^3))$.
Furthermore this $T$ depends only on $\|u_0\|_{H^1}$.
\end{prop}
\begin{proof}
We prove this using a standard contraction argument. We adopt the notation $\| \cdot \| = \| \cdot\|_{L^\infty([0,T]H^1(\R^3))}$:
\begin{align*}
\|F(u)& - F(v)\| \le \left\|\int_0^t e^{i(t-s)\Delta} \left[N(u(s)) -
N(v(s))\right] ds\right\| + \left\|\int_0^t e^{i(t-s)\Delta} \left[Vu(s) - Vv(s)\right] ds\right\|\\
& \le c\left(\|N(u(t)) - N(v(t))\|_{L^2_{[0,T]}W^{1,6/5}_x} + T\|Vu(t) - Vv(t)\|\right)
\end{align*}
But then
\begin{align*}
c\|N(u(&t)) - N(v(t))\|_{L^2_{[0,T]}W^{1,6/5}_x} \le c T^{1/2}\|N(u)-N(v)\|_{L^\infty_{[0,T]}W^{1,6/5}_x}\\
&\le c T^{1/2}\Big[\left\|\left(|x|^{-1}*|u|^2\right)(u -
v)\right\|_{L^\infty_{[0,T]}W^{1,6/5}_x} + \left\|\left(|x|^{-1}*u(\bar u -
\bar v)\right)v\right\|_{L^\infty_{[0,T]}W^{1,6/5}_x} +
 \\
&\qquad\qquad\qquad\qquad\qquad\qquad\qquad\qquad\qquad\qquad\qquad\left\|\left(|x|^{-1}*(u-v)\bar
v\right)v\right\|_{L^\infty_{[0,T]}W^{1,6/5}_x}\Big]\\
 &\le c T^{1/2}\|u-v\|\left(\|u\|^2 + \|u\|\|v\| +
\|v\|^2\right)
\end{align*}
Thus taking
\[T^{1/2} \le \frac{1}{c\left(\|u\|^2 + \|u\|\|v\| + \|v\|^2 + \|V\|_{W^{1,\infty}(\R^3)}\right)},\]
we find that $F$ is a contraction on a
closed ball of $L^\infty([0,T],H^1(\R^3))$, implying there exists a
solution to \eqref{e:hartree}.
\end{proof}

We then use almost conservation of energy to extend this to global well-posedness.

\begin{prop}
The equation \eqref{e:hartree} has a solution in $L^\infty(\R,H^1(\R^3))$
\end{prop}

\begin{proof}
Because of Proposition \ref{p:locwp}, it is sufficient to prove that
the $H^1$ norm of $u$ is bounded. Clearly $\|u\|_{L^2}$ is preserved
so it suffice to bound $\| \nabla u\|_{L^2}$. To do this we study
the energy
\[
E(t)=\|\nabla u\| - \int_{\R^3} N(u)\overline{u}.
\]

An argument as above shows that

\[
\int \left(|x|^{-1} * |u|^2 \right) |u|^2  \le \||x|^{-1} * |u|^2\|_{L^3} \|u^2\|_{L^{3/2}}\le c\|u\|_{L^2}^3 \|\nabla u \|_{L^2} \le \frac{c}{\epsilon} \|u\|_{L^2}^3 + c\epsilon\|\nabla u \|_{L^2}.\]

From this we deduce that
\[
\| \nabla u\|^2_{L^2} \le c \left( E(t) +
\|u\|_{L^2}^3+\|V\|_{W^{1,\infty}} \right).
\]
This bounds $\|u\|_{H^1_x}$ uniformly in time, giving the desired conclusion.\end{proof}

\textbf{Acknowledgments.} We would like to thank Maciej Zworski for suggesting this problem and for his generous guidance and advice. Thanks also to Justin Holmer and Sebastian Herr for helpful discussions regarding Appendix B. Finally, we are grateful for support from NSF grant DMS-0654436.

\end{document}